\documentclass[11pt,a4paper,draft]{article}
\usepackage{amsmath,mathrsfs,amssymb,amsthm,amscd}
\usepackage[]{fontenc}
\usepackage{xy}
\usepackage{enumerate}
\xyoption{all}

\numberwithin{equation}{section}

\theoremstyle{plain}
\newtheorem{theorem}[equation]{Theorem}
\newtheorem{corollary}[equation]{Corollary}
\newtheorem{proposition}[equation]{Proposition}
\newtheorem{lemma}[equation]{Lemma}

\theoremstyle{definition}
\newtheorem{definition}[equation]{Definition}

\newtheorem{remark}[equation]{Remark}

\begin{document}

\title{Bounds of the number of rational maps\\
between varieties of general type
}

\author{
\addtocounter{footnote}{5} J.C. Naranjo \footnote{Partially
supported by the Proyecto de Investigaci\'on BFM2003-02914 and
ACI-2003} ,
 G.P. Pirola
\footnote{Partially supported by
1) PRIN 2003 {\em "Spazi di moduli e teoria di Lie"}; 2) Gnsaga;
3) Far 2004 (PV):
 {\em Variet\`{a} algebriche, calcolo algebrico, grafi orientati e topologici}.}
}

\maketitle

\pagestyle{plain} \setcounter{tocdepth}{1}

\begin{abstract}{We give a bound for the number of  rational maps between
algebraic varieties of general type under mild hypothesis on the
canonical map. We use an idea inspired by Tanabe's work. Instead
of attaching a morphism of Hodge structures to a rational map  we
 simply associate to it  a piece of the integral Hodge lattice. This procedure does
not give an injective map, but by means of a geometric argument,
we can estimate the number of maps with the same image. \par
}
\end{abstract}

\section{Introduction}

De Franchis  proved in 1913 (see \cite{dF}) that the set of morphisms between  two Riemann
surfaces of genus at least $2$ is finite. In other words, he showed the finiteness of the set
$$
M(X,Z)=\{f:X \longrightarrow Z\,|\, f \text{ non-constant}\},
$$
where $X,Z$ are curves of genus at least $2$.

Martens (cf. \cite{Ma}) gave an effective bound of the number of elements $m(X,Z)$ of
this set. Other estimates can be deduced from the effective bounds
for the number elements of
$$
M(X)=\{f:X\longrightarrow Y\,|\, f \text{ non-constant},\, Y
\text{smooth curve of genus}\ge 2\},
$$
obtained in \cite{HS:dF},\cite{Ka} and \cite{AP:dF}.

Probably the most interesting open problem in the topic (see
\cite{He}) is whether $m(X,Z)$ can be bounded by a polynomial on
the genus of the curves. Kani (\cite{Ka}) showed that this is not
true for $M(X)$.

In \cite{Ta} Tanabe has improved the known bounds for $m(X,Z)$. In
all the previous proofs morphisms of homological lattices were
used to represent maps, or even correspondences, on curves. The
main idea of Tanabe's work is to represent each map by a single
element of the singular homology group of $X$. This enables him to
control, with a geometric argument, the number of maps represented
by the same element. In fact, his proof can be separated into two
parts. In the first he shows that, fixing a holomorphic form
$\alpha $ on the target, two maps in which the pull-backs of
$\alpha $ are the same ``differ" in a finite number of choices
depending polynomially on the genus.

In the second part he assumes that $\alpha $ is the $(1,0)$-part
of an element $\tilde {\alpha}$ of the integral lattice with
minimal norm. Then he attaches to the map the  pull-back of
$\tilde {\alpha }$. This is an element in the integral lattice of
the source. He then shows that we can reduce the lattice modulo
$d$ with $d$ greater than twice the degree of the map, without
losing information.

We refer to the first
argument as the ``geometric part" and to the second as the
``torsion part" of Tanabe's work.

This paper concerns the same problem in a higher dimension, that
is, we consider
$$
M(X,Z)=\{f:X - - \to Z\,|\, f \text{ rational dominant map }\},
$$
for $X,Z$ varieties of general type of the same dimension. As
above, $m(X,Z)$ denotes the number of elements of $M(X,Z).$

It was proved by Kobayashi-Ochiai that $m(X,Z)$ is finite (see
\cite{KO} and also \cite{BD}). Moreover there is an effective bound  in \cite{He} for
complex manifolds with ample canonical bundle obtained by means of
Chow varieties. This method provides necessarily a bound with a
very high exponential.

In this paper we use an idea inspired by Tanabe's work. Instead of
attaching a morphism of Hodge structures to a rational map we
simply associate to it a piece of the integral Hodge lattice. This
procedure does not give an injective map, but, by means of a
geometric argument  we can estimate  the number of maps with the
same image.

We do not need the restrictive hypothesis which guarantees the
injectivity of the representation of the elements of $M(X,Z)$ as
maps of Hodge structures. We can thus  find good bounds under weak
hypotheses. In fact, we find much better bounds for
$n$-dimensional varieties than the ones currently known.

We use two approaches. The first works in dimensions $1$, $2$ and
$3$ and gives better results. The second applies in any dimension,
under a more restrictive hypothesis.

Now we explain the ideas of the proofs:
 First we generalize the geometric part of Tanabe's work to
surfaces with $p_g$ at least $2$ by using appropriate pencils of
$2$-forms on $Z$. Since $m(X,Z)$ is a birational invariant we may
assume that $X$ and $Z$ are minimal. Next we represent the map
using couples of elements in the transcendental lattice of the
source variety. Roughly speaking, the transcendental lattice is
the complement of the Neron-Severi group in the second cohomology
group of the surface. The geometric part allows us to estimate the
number of maps which are represented by the same couple of
elements of the lattice. To do this we use the following fact,
which is elementary but very useful: there exists an open set where all the
maps are well-defined and such that for each point of this open
set two different maps take different values. Then, by using the
fact that the curves are moving in a pencil and thus cut this open
set, one can reduce the proof to the one-dimensional case.

Next, instead of the torsion lemma we use a packing lemma due to
Kani. To do so we give a lower estimate for the distance of two
different elements. We obtain a bound of $m(X,Z)$ in terms of
$K_X^2$, $K_Z^2$ and the second Betti number $b_2(X)$ (see
\ref{th2}). By combining Bogomolov-Miyaoka-Yau and Noether
inequalities one can obtain an estimate in terms of the Euler
characteristic (see \ref{cor2}).

Observe that, since we are not assuming that $X,Z$ are canonical,
the representation of the maps in $M(X,Z)$
as maps of transcendental Hodge structures is not injective in general.

Note also that using the packing lemma (instead of the torsion
lemma) in the 1-dimensional case, we obtain a result which is
slightly better than  Tanabe's (see \ref{th1}).

Then,  with some additional hypotheses we can give a bound for
threefolds following a similar argument. Note the difference in
the arguments for surfaces and threefolds. In the first case, to
prove the geometric lemma we reduce the proof to the one for
curves. Instead, in the case of threefolds
 we need to use the full
result on surfaces.

Apparently this ``inductive procedure" does not  extend to higher
dimensions due to the method used and to the lack of a smooth
minimal model in higher dimension.

In the last section we extend the torsion  part in Tanabe's work.
We use this to give a bound in general (see \ref{thn}). This bound
is clearly worse than the one obtained  for surfaces and
threefolds.

The paper is organized as follows: in \S 2 we give some
preliminaries, mainly on Hodge structures. We also recall Kani's
packing lemma.

To give the statements of the following theorems, we introduce the
following function
$$
P(a,e)=(a+1)^e-(a-1)^e, \qquad a\in \mathbb R,\, e\in \mathbb N.
$$
This is a polynomial on $a$. Its leading term is $$2ea^{e-1}.$$

We also denote
$$
\rho =\rho (X,Z) =\frac {K^n_X}{K_Z^n}
$$
where $X,Z$ are $n$-dimensional varieties (if $n=1$, then $\rho=\frac{g(X)-1}{g(Z)-1}$).

Let $b_i(X)$ be the Betti number $\dim H^i(X,\mathbb C)$.

In \S 3 we prove the
theorem on curves:

\vskip 5mm
\begin{theorem}\label{th1}
Let $X,Z$ be smooth irreducible projective curves of genus $\ge
2$. Then
$$
\begin{aligned}
m(X,Z)&\le 4(g(X)-1)\rho\,P(2\rho,2g(X))=
\\
= &8(g(X)-1)\,\rho\,[ \binom{2g(X)}{1} (2\rho )^{2g(X)-1}
+\binom{2g(X)}{3} (2\rho )^{2g(X)-3} \dots
 ].
 \end{aligned}
$$
\end{theorem}

Section 4 is devoted to the proof of:

\vskip 5mm
\begin{theorem}\label{th2}
Let $X,Z$ be smooth irreducible projective minimal surfaces of
general type. Assume $p_g(Z) \ge 2$. Then
$$
m(X,Z)\le 4 (K_X^2)^2\, P(4\sqrt 2 \,\rho ,2b_2(X)-2).
$$
\end{theorem}

\vskip 5mm
 Since $\rho \le K_X^2\le 9\chi (\mathcal O_X)$
(Bogomolov-Miyaoka-Yau) and
$$
\begin{aligned}
b_2(X)&=\chi _{\text{top}}(X)+2q(X)-2 =12 \chi (\mathcal O_X)-K_X^2
+2q(X)-2=
\\
&=10 \chi (\mathcal O_X)+2p_g(X)-K_X^2\le 10 \chi (\mathcal O_X) + 4
\end{aligned}
$$
(we use Noether's formula and Noether's inequality) we immediately
obtain a bound for surfaces in terms of the Euler characteristic:

\vskip 5mm
\begin{corollary}\label{cor2}
Let $X,Z$ be smooth irreducible projective surfaces of
general type. Assume $p_g(Z) \ge 2$. Put $\chi  =\chi (\mathcal O_X)$. Then
$$
m(X,Z)\le 4\cdot 9^2\,\chi ^2\,P(36\sqrt 2 \,\chi ,20\chi
+6).
$$
\end{corollary}

In \ref{cor2} we do not assume the surfaces are minimal. This will
be useful in the proof of the next theorem for threefolds, which
will be given in section 5.
 \vskip 5mm
\begin{theorem}\label{th3}
Let $X,Z$ be smooth irreducible projective complex threefolds of
general type. Assume that $K_X,\,K_Z$ are nef, $p_g(Z)\ge 2$ and the
image of $Z$ by the bicanonical map has dimension at least $2$. Then
$$
m(X,Z)\le  4\cdot 9^2\, h^2\, K^3_X\, P(36\sqrt 2 \,h,20 h +6)\cdot
        P(4\sqrt 2 \rho,2b_3(X)),
$$
where
$h=h^0(X,\mathcal O(2K_X))+h^0(X,\Omega ^2_X)-p_g(X)+1$.
\end{theorem}

\vskip 5mm
 Finally in \S 6 we find a bound for $n$-dimensional varieties by extending
 Tanabe's torsion part to higher dimension.
The result we obtain is:

\vskip 5mm
\begin{theorem}\label{thn}
Let $X,Z$ be two $n$-dimensional varieties of general type such that  $K_X,\,K_Z$
are nef and $\dim (\varphi _{|K_Z|}(Z))\ge n-1$. Then:
$$
m(X,Z)\le 2n (K_X^n)^2(2\rho +1)^{b_n(Z)\cdot b_n(X)}.
$$
\end{theorem}

In the case of birational automorphisms we obtain a bound with a
lower exponent than the one given in \cite{He},

\vskip 5mm
\begin{corollary}
Let $X$ be a variety of general type with $K_X$ nef and such that
$\dim (\varphi _{|K_X|}(X))\ge n-1$, then
$$
\# \text{aut} (X)\le 2n \,(K_X^n)^2\, 3^{b_n(X)^2}.
$$
\end{corollary}

\vskip 5mm \textbf{Acknowledgements.} We are grateful to M.A.Barja and A.Collino for
valuable suggestions during the preparation of this paper.

The work was completed during the second author's stay at the
Institut de Matem\`{a}tiques de la Universitat de Barcelona
(IMUB). This stay was supported by the Generalitat de Catalunya's
PIV programm.

\vskip 10mm
\section{Notations and Preliminaries}

\subsection{Notations}

Throughout the paper we use the symbols $M(X,Z)$ and $m(X,Z)$ as
in the introduction: the former is the set of rational dominant
maps from to $X$ to $Z$ and the latter is its number of elements.

Analogously, $M_r(X,Z)$ is the subset of maps with fixed degree
$r$ and the number of its elements is $m_r(X,Z)$.

We also keep the notations:
$$
\rho=\rho (X,Z)=\frac {K^n_X}{K_Z^n}, \qquad \text{where }\, n=\dim X=\dim Z,
$$
and
$$
b_i(X)=\dim H^i (X,\mathbb C).
$$

We work over the complex numbers.
In this paper variety means irreducible, smooth, projective, complex variety.

\subsection{Hermitian spaces }
Let $(V,h_{V})$ be a hermitian of finite dimension space
 we shall denote the norms for $v\in V$
 by
$$\|v\|=\sqrt{h_{V}(v,v)}.$$
We recall that, for any map $f:V\to W$ between twohermitian
spaces, we may define the adjoint map $g: W \to V$ as the unique
 linear map such for all $v\in V$ and $w\in W$
 $$
 h_{W}(f(v),w)=h_{V}(v,g(w)).
 $$

\vskip 3mm
\subsection{Hodge structures and morphisms}

Let $X$ be a complete smooth algebraic variety of dimension $n.$
Let $H^{n}(X)$ be the Hodge structure on $X$ on the $n$-cohomology
of $X.$ The lattice $H_{\mathbb Z},$ is
$$H^{n}(X,\mathbb Z)/torsion$$ and the standard Hodge
decomposition : $H_{\mathbb Z} \otimes \mathbb C= H^{n}(X,\mathbb
C)=\oplus H^{p,q},$ $ H^{p,q} =\overline  {H^{q,p}}.$ Integration
gives a natural polarization:
$$
Q: H^{n}(X,\mathbb C)\times H^{n}(X,\mathbb C) \to \mathbb C
$$
which is unimodular, by Poincar\'e duality, on $H_{\mathbb Z}.$
We recall that a Hodge substructure $R$ of $H$ is given by a sublattice
$R_{\mathbb Z}\subset H_{\mathbb Z} $ such that
$R_{\mathbb C}=R_{\mathbb Z}\otimes\mathbb C=\oplus R^{p,q},$  where
$ R^{p,q}= H^{p,q}\cap R_{\mathbb C}.$ The restriction of $Q$
gives a polarization of  $R$ non necessarily unimodular over the integers.
The polarization makes possible to define the orthogonal Hodge structure
$R'.$ Set
$R'_{\mathbb Z}=\{\gamma \in H_{\mathbb Z}\,|\, Q(\gamma, \beta)=0\
  \forall \beta\in R_{\mathbb Z}\}.$
One has $R_{\mathbb C}\oplus R'_{\mathbb C}= H_{\mathbb C}.$

\begin{definition}\label {birhodge}
The {\bf transcendental} Hodge structure of $X$ is the smallest
Hodge substructure $T_X$ of $H^n(X)$ containing $H^{n,0}(X)$. Its
lattice $T_{\mathbb Z,X}$ will be called the transcendental
lattice of $X.$ For any integer $d\geq 2$ let $T_{d,X}=T_{\mathbb
Z, X}/d\cdot T_{\mathbb Z,X}.$ Observe that if $n$ is even then
there exists a $(\frac n2,\frac n2)$-integral class induced by a
projective immersion; therefore $T_X\ne H^n(X,\mathbb C)$ and
$\dim T_X \le b_n(X)-1$.
\end{definition}
Note that $T_{X}$ is a birational invariant of $X.$ Since $T_X$ is
contained in the primitive cohomology, then  (due to the
Hodge-Riemann relations, see \cite{GH}, page 123) the cup-product
modified with the Weil operator induces on $T_X$ a hermitian
product that we denote simply by $(,)$.

\vskip 3mm
Let $Z$ be another smooth complete variety of dimension $n$ and
$$f:X --\to Z$$
be a dominant rational map of degree $r=\deg f$. We then have two
Hodge structure morphisms:
$$
f^{\ast}:T_{Z}\to T_{X},\quad
f_{\ast}:T_{X}\to T_{Z}
$$
We have that they are adjoint maps; in other words:
\begin{equation}
(\alpha, f_{\ast}(\beta))=(f^{\ast}(\alpha),\beta).
\end{equation}

We may also define the group
 homomorphism $f_{d}: T_{d,Z}\to T_{d,X}$ induced by $f^{\ast}$.

We have the following:
\begin{lemma}\label{norma}\
\begin{enumerate}
\item [a)] If $\gamma\in T_{Z}$ then $f_{\ast}f^{\ast}(\gamma)=r
\gamma$.
\item [b)] If $\gamma\in T_{Z}$ then $ \| f^{\ast}(\gamma)\|=\sqrt {r}\,\|\gamma\|$
\item [c)] For any $ \beta\in T_{X} $
$\| f_{\ast}(\beta)\|\le \sqrt r \,\|\beta\|$
and $\| f_{\ast}(\beta)\|= \sqrt r \,\|\beta\|$ if and only if
there is $\gamma\in T_{Z}: \beta=f^{\ast}(\gamma).$
\end{enumerate}
\end{lemma}

\begin{proof}
Parts a) and b) are well known. To see c), we
write $\beta =f^{\ast }(\beta _0)+\eta$, where $\eta $ is orthogonal to the image of $f^{\ast }$. Therefore:
$$
\begin{aligned}
&(f_{\ast }(f^{\ast }(\beta _0)+\eta),
 f_{\ast }(f^{\ast }(\beta _0)+\eta))=
(f_{\ast }(f^{\ast }(\beta _0)),
 f_{\ast }(f^{\ast }(\beta _0)))=
\\
&deg(f)^2(\beta _0,\beta _0)=deg(f)\,(f^{\ast}(\beta_0),
f^{\ast}(\beta_0))\le deg(f)\,(\beta,\beta).
\end{aligned}
$$
\end{proof}

\subsection{Packing lemma}

We will need the following lemma, which appears in \cite{Ka}. To
state it more clearly we define:
$$
P(a,e)=(a+1)^e-(a-1)^e,
$$
where $a\in \mathbb R$ and $e\in \mathbb N$. Observe that $a\le
a'$ implies $P(a,e)\le P(a',e)$. Also $e\le e'$ implies $P(a,e)\le
P(a,e')$.

\begin{lemma}\label{kannipl}
Let $v_1,\dots,v_N\in \mathbb R^v$, $\|v_i\|=R>0,\, \forall i$.
Assume $\|v_i-v_j\|\ge 2d,\, \forall i,j, i\ne j$, then
$$
N\le P(\frac Rd,v)=2[\binom{v}{1}
(\frac Rd)^{v-1}+
\binom{v}{3} (\frac Rd)^{v-3}+
\dots
 ]
$$
\end{lemma}

\subsection{Degree of rational maps}
 Let $X,Z$ be two $n$-dimensional varieties of general type such that $K_X$ and $K_Z$ are nef. One has

\begin{lemma} \label{degrat}
Let $f:X - - \to Z $ be a rational dominant map. Then
$$
\text{deg }(f)\le \rho (X,Z).
$$
\end{lemma}
\begin{proof}
For $n=1$ it is a consequence of Riemann-Hurwitz formula. Assume
$n\ge 2$. Since $K_X,\,K_Z$ are nef, by taking $l>>0$,
 the linear systems  $|lK_Z|$, $|lK_X|$ are base-point-free. Then, we can think of $f$ as a
linear projection in a projective space. Then the degree of $f$ is
bounded by the quotient of the degrees of $\varphi _{|lK_X|}(X)$
and $\varphi _{|lK_Z|}(Z)$; hence $deg(f)\le K_X^n/K_Z^n$.
\end{proof}

\subsection{Rational domain}

With two $n$-dimensional varieties of general type $X,Z$ fixed,
recall that $M(X,Z)$ is finite (see \cite{KO}). Assuming that it
is not empty, we label its elements $ M(X,Z)= \{f_{i}\},$ \
$i=1,\dots,m(X,Z).$

\begin{definition} A Zariski open set of $W\subset X$ will be called
a rational domain for $X$ and $Z$ if  any $f_{i}\in M(X,Z)$
defines a regular morphism $f_{i|W}:W\to Z$ and for any point
$x\in W$ $f_{i}(x)=f_{j}(x)$ implies $i=j.$ \label{ratdom}
\end{definition}

A rational domain always exists  since
 the closure of the sets
$f_{i}=f_{j},$ $i\neq j$ are proper algebraic subsets of $X$. Note
for $x\in W$,
$$\#\{z_{i}=f_{i}(x)\}=m(X,Z).$$

\vskip 10mm
\section{Curves}

We consider the case of curves, so  $ 1=\dim X =\dim Z$.

\subsection{Tanabe's geometric lemma}

Our first goal is to rewrite the geometric part of Tanabe (see \cite{Ta}). We fix a holomorphic form on $Z$,
 $ 0\neq \alpha \in H^0(Z,K_{Z})$ and we say that two maps $f,g$ are equivalent if and only if
$f^{\ast }(\alpha )=g^{\ast }(\alpha )$.
We want to  give a bound of the number elements of the equivalence class $[f]$ of a map $f$.

Let $x$ be a zero of $f^{\ast }(\alpha )$ and put $z=f(x)\in Z$. Let us denote by $\mathbb D$  the Poincar\'e
 disk and $p:\mathbb D\to X$ and $q:\mathbb D\to Z$ be the
 universal coverings such that $ p(0)=x$ and $q(0)=z.$ To any holomorphic map $f:X\to Z$
such that $f(x)=z,$ there is a unique
lifting holomorphic map $F:\mathbb  D\to \mathbb D$ such $F(0)=0$ and
$q(F(t))=f(p(t)) $ for all $t\in \mathbb D.$
Assume $g\in [f]$  is another non-constant holomorphic function with
 $g(x)=z$ and lifting $G:\mathbb D\to \mathbb D$, $G(0)=0.$

We give the following global version of Tanabe's lemma.
\begin{lemma}
Under the previous hypothesis,
\begin{enumerate}
\item [a)] There is a constant $c$ such that $F(t)=c\, G(t).$
\item [b)] If $n$ is the order of $ \alpha $ at $z$ then $c^{n+1}=
1.$
\end{enumerate}
\end{lemma}
\begin{proof}
Let us consider
 the pull-back of the form $\alpha $ on
$\mathbb D$:
$$
q^{\ast}(\alpha)= k(t)dt.
$$
The condition $f^{\ast}(\alpha ) = g^{\ast}(\alpha )$
 translates into
$$ k(F(t))\cdot dF(t)= k(G(t))\cdot dG(t).$$
If $K(t): D\to \mathbb R$ is the primitive of $k(t)$ such that
$K(0)=0$ we obtain:
$$ K(F(t))= K(G(t)).$$
Now if $k(t)$ has  order $n$ at zero, $K(t)$ has a zero of order
$n+1$ and we can find a function $w(t)$ defined near zero such
that $w(K(t))=t^{n+1}.$ From $w(K(F(t)))= w(K(G(t)))$ we obtain
$$
F(t)^{n+1}=G(t)^{n+1}.
$$
That is, $F(t)=c\, G(t)$, $c^{n+1} = 1.$
 \end{proof}

\begin{corollary}\label{tgl1}
The number of elements of $[f]$
is less than or equal to $4(g(X)-1)$.
\end{corollary}

\begin{proof}
Due to the lemma, for each  zero $x$ of $f^{\ast }(\alpha )$ we have at most $$ord_{f(x)}(\alpha )+1$$
maps of $[f]$ with the same image at $x$. Consider for any $x\in (f^{\ast }(\alpha))_0$ the set
$$
A_x\,=\,\{ g\in [f]\,|\, g(x)=z  \},
$$
where $z$ is a fixed zero of $\alpha $. Observe that
$$
[f]\,=\,\bigcup _{x\in (f^{\ast }(\alpha))_0}A_x.
$$
Therefore
$$
\begin{aligned}
\#[f] &\le \sum _{x\in (f^{\ast }(\alpha))_0}(ord_{f(x)}(\alpha )+1)
\le \sum _{x\in (f^{\ast }(\alpha))_0}(ord_x(f^{\ast }(\alpha ))+1)\\
&\le 2g(X)-2+\sum _{x\in (f^{\ast }(\alpha))_0}1\le 4(g(X)-1).
\end{aligned}
$$
\end{proof}

\subsection{Proof of Theorem \ref{th1}}

Let $\tilde{\alpha}$ be a non-trivial element in $T_{Z,\mathbb Z}$
with minimal norm. We denote by $\alpha $ its $(1,0)$-part. So
$$
\tilde{\alpha}=\alpha +\overline{\alpha }.
$$
We define the equivalence relation $\sim $ in $M(X,Z)$ as follows:
$$
f\sim g\,\text{ if and only if }\, f^{\ast }(\tilde {\alpha })=
g^{\ast }(\tilde {\alpha }).
$$

It is obvious that
$$
f\sim g\,\text{ if and only if }\, f^{\ast }(\alpha )=
g^{\ast }( \alpha ).
$$
In particular, the class of $f$ under the relation $\sim $ is the set $[f]$ considered in the section 3.1.

Let us fix a positive integer $r$. Observe that $\sim$ is in fact an equivalence relation in $M_r(X,Z)$,
since $\|f^{\ast}(\tilde {\alpha })\|=\sqrt {deg(f)}\,\|\tilde {\alpha }\|$. By \ref{degrat} the constant $r$
is bounded above by $
\rho $. So, due to \ref{tgl1}, we get

$$
m(X,Z)=\sum _{r=1}^{\rho}m_r(X,Z)\le 4(g(X)-1)
\sum _{r=1}^{\rho }\#(M_r(X,Z)/\sim).
$$
Now we use the injection
$$
\begin{aligned}
M_r(X,Z)/\sim &\hookrightarrow H^1(X,\mathbb Z) \otimes \mathbb R\\
f&\longmapsto v_f:=(1/\|\tilde {\alpha }\|)\,f^{\ast}(\tilde {\alpha })
\end{aligned}
$$
to bound the number of elements of the quotient $M_r(X,Z)/\sim$. Observe
that the image belongs to the sphere of radius $\sqrt {r}$
centered at the origin in a real vector space of dimension
$2g(X)$.

\begin{proposition}\label{dist}
Let $f,g:X\longrightarrow Z$ be two maps of degree $r$ such that
$f^{\ast }(\alpha )\ne g^{\ast}(\alpha )$. Then
$$\|v_f-v_g\|\ge \frac 1{\sqrt r}. $$
\end{proposition}

\begin{proof}
Observe that
$$
((f_{\ast }-g_{\ast})(f^{\ast}(\tilde {\alpha })-g^{\ast}(\tilde {\alpha } )),\tilde {\alpha } )=
 (f^{\ast}(\tilde {\alpha } )-g^{\ast}(\tilde {\alpha } ),f^{\ast}(\tilde {\alpha })-g^{\ast}(\tilde {\alpha } ))
 \ne 0,
$$
hence we can assume $f_{\ast }(f^{\ast}(\tilde {\alpha })-g^{\ast}(\tilde {\alpha }))\ne 0.$
Therefore by using the minimality of the norm of $\tilde {\alpha }$:
$$
\|\tilde {\alpha } \|\le \|f_{\ast }(f^{\ast}(\tilde {\alpha })-g^{\ast}(\tilde {\alpha }))
\|\le \sqrt r\,\| f^{\ast}(\tilde {\alpha })-g^{\ast}(\tilde {\alpha })  \|,
$$
which implies the statement.
\end{proof}

By Lemma \ref{kannipl} with $d=\frac 1{2\sqrt r}$, $v=2g(X)$ and $R=\sqrt r$, we get
$$
\#(M_r(X,Z)/\sim) \le P(2r,2g(X)) \le P(2\rho , 2g(X)).
$$

Together, this gives
$$
m(X,Z)\le 4(g(X)-1)\,\rho \,P(2\rho, 2g(X)),
$$
proving \ref{th1}.

\begin{remark}
Notice that
$$
P(2\rho,2g(X))=[\binom{2g(X)}{1} (2\rho )^{2g(X)-1}+
\binom{g(X)}{3} (2\rho)^{2g(X)-3} +\dots ],
$$
so the dominant term of the bound has exponent $2g(X)-1$ instead
of the exponent $2g(X)$ that appears in Tanabe's bound.
\end{remark}

\begin{remark}
One can improve the bound given above by finding  a better lower bound for $\|v_f-v_g\|$. In fact we can prove:
 $\|v_f-v_g\|\ge \sqrt {\frac{r^2+1}{r^3}}$.
\end{remark}

\vskip 10mm
\section{Surfaces}

In this section we assume that $X$ and $Z$ are surfaces of general type and that $p_g(Z)\ge 2$.
The general strategy to find a bound
for $m(X,Z)$ is similar to that used for curves: we find a bound
for the number of maps which fix a pencil of $2$-holomorphic forms
minimal in some sense. Then we use the transcendental lattice to
represent the maps and we prove a result similar to
\ref{dist}.

\subsection{Generalization of the geometric lemma}

We fix two independent $(2,0)$ forms $\alpha$ and $\beta$ on $Z.$
We define the following equivalence relation on  $M(Z,X):$
$$ \label{sim}
 f\sim g\  \iff f^{\ast}(\alpha)= g^{\ast}(\alpha) \text { and }
f^{\ast}(\beta)= g^{\ast}(\beta).
$$

\begin{remark} \label{eqind}
 If $f\sim g$ then $ |f^{\ast}(\beta)|^{2}= |g^{\ast}(\beta)|^{2}$
 then $\deg f=\deg g$, that is, the above relation gives a
 equivalence relation on $M_{r}(X,Z).$
 \end{remark}

As in section 3, we would like to evaluate the number of elements
 in an equivalence class $[f]$.
To do so we take the pencil $L$ generated by $\alpha$ and $\beta.$
We also let $B$ be the base curve (it could be $B=\emptyset$) of
$L$. We may assume that $\beta$ is the general element of $L$.
 Then the zero divisor
 $(\beta)_0$ of $\beta$ can be written
 as
 $$
 (\beta)_0= B+\sum_{i=1}^{s}  C_{i},
 $$
 where $C_{i}$ are reduced and irreducible
 curve of geometric genus
 $g$ with $g\geq 2.$ We may also assume that
 $C_{i}\cdot C_{j}\geq 0$.

 Now we denote by $L'$ the pencil
 $f^{\ast}(L)$, which is independent of the choice of a map in $[f]$.

Then we obtain
$$
  (f^{\ast }(\beta ))_0= B'+\sum_{i=1}^{s'}  C'_{i},
$$
 where $B'$ is the base divisor of the pencil and $C'_i$ are irreducible reduced curves of genus $g'\ge 2$.
 We denote by $N_i$ (respectively $N_i'$) the normalization of the curve $C_i$ (respectively $C_i'$).
 We have the following lemma

\begin{lemma} \label{genN} Let  $s$  be
  the number  of irreducible components of
  $(\beta)_0- B$.
  Then:
  \begin{enumerate}
 \item[a)] $s\leq K^{2}_{Z}$.
 \item[b)] $g(N_i)\leq K^{2}_{Z}+1$, $g(N_i')\leq K^{2}_{X}+1$.
\end{enumerate}
  \end{lemma}
\begin{proof}
  \begin{enumerate}
  \item[a)]
Since $C_{i}$ moves, then $K_Z\cdot C_i\ge 1$. Therefore:
$$
s\leq K_{Z}\cdot \sum_{i=1}^{s}  C_{i}=K_Z\cdot (K_Z-B)\le
K^{2}_{Z}.
$$
 \item[b)]
 The proof is given on $Z$. Observe that,
since $K_Z$ is nef:
$$
(K_Z+C_i)(K_Z-C_i)\ge C_i(K_Z-C_i)\ge 0.
$$
So, $K_Z^2\ge C_i^2$.
In fact, if there is more than one component, by $2$-connectivity
$C_i(K_Z-C_i)\ge 2$ and then $K_Z^2\ge C_i^2+2$.

 Then we have
$$
g(N_i)\le p_a(C_i)=\frac 12 (C_i^2+C_i\cdot K_Z)+1\le \frac
12(K_Z^2+K_Z^2)+1=K_Z^2+1.
$$
\end{enumerate}
\end{proof}

Let us consider $Z'\longrightarrow \mathbb P^1$ to be the minimal
resolution of the pencil $$Z- - \to \mathbb P^1.$$ Let $X'- -
\to\mathbb P^1$ be the minimal resolution of the pencil on
$X\times _Z Z'$. Then the map $f$ and the forms $\alpha $, $\beta
$ pull-back to $f'$, $\alpha '$, $\beta '$  and we obtain
$$[f]=[f']=\{g':X' - - \to Z' \,|\, f'^{\ast }(\alpha ')=g'^{\ast
}(\alpha '),f'^{\ast }(\beta ')=g'^{\ast }(\beta ')\}.$$

Observe that an irreducible component of a general fibre of the
pencil on $X'$ (resp.$Z'$) is $N'_i$ (resp. $N_i$).

Now we fix the component $N_1'$.
 Then $[f']$ is the union of the subsets of maps which send $N'_{1}$ to $N_{i}$, $i=1,\dots ,s$:
$$
[f']= \bigcup _i\{g\in [f'] \, |\, g(N_1')=N_i\}.
$$
Observe that $N_1'$ intersects the rational domain of $X'$ and
$Z'$ (see \ref{ratdom}) because it is a component of a generic
element of a pencil. Moreover by taking the residue of $\alpha
'\otimes \alpha '/\beta'$ along $N_i$ a $1$-form $\hat {\alpha } _i$ is
induced on $N_i$ (see \cite {GH},pp.500-505). By definition, the
pull-back of $\hat {\alpha } _i$ is the same for all the maps in
$[f]$. Therefore
$$
\{g\in [f'] \, |\, g(N_1')=N_i\} \subset
\{g:N_1'\longrightarrow N_i\,|\, g^{\ast }(\hat {\alpha } _i)
\text{ fixed}\}.
$$

We are ready to prove

\begin{proposition}\label{tgl2}
One has the inequality:
$$
\#[f]\le 4K_Z^2\,K_X^2.
$$
\end{proposition}
\begin{proof}
We use \ref{tgl1} in the last inclusion of sets and we obtain, by means of \ref{genN}:
$$
\#[f]=\#[f']\le\sum_1^s 4(g(N'_1)-1) =  4 s\,K_X^2\le 4
K_Z^2\,K_X^2.
$$
\end{proof}

\subsection{Proof of \ref{th2}}

Let $\tilde {\alpha}\in T_{Z,\mathbb Z}$ be an element of the
transcendental lattice in $Z$ with minimal norm (see
\ref{birhodge}). Put $\alpha $ the $(2,0)$-component of
$\tilde{\alpha}$. The smallest Hodge substructure containing
$\tilde{\alpha}$ is denoted by $\left\langle
\tilde{\alpha}\right\rangle $. If $\left\langle
\tilde{\alpha}\right\rangle =T_Z$, then $\tilde {\beta }$ is any
$(2,0)$-form linearly independent with $\alpha$. If, on the
contrary, $\left\langle \tilde{\alpha}\right\rangle \ne T_Z$ we
can find a decomposition of Hodge structures $T_Z=\left\langle
\tilde{\alpha}\right\rangle \oplus ^{\perp  }R.$ Then we choose an
element $\tilde {\beta }\in R_{\mathbb Z}$ with minimal norm. By
construction its $(2,0)$-component $\beta $ is linearly
independent with $\alpha$.

\begin{definition}
Two rational maps $f,g\in M(X,Z)$ are equivalent if and only if
$f^{\ast }(\tilde {\alpha})=g^{\ast }(\tilde {\alpha})$
and $f^{\ast }(\tilde {\beta})=g^{\ast }(\tilde {\beta }).$
\end{definition}

We denote this relation also by $\sim $, since by the next lemma
it coincides with the relation given in \ref{sim}.

\begin{lemma}\label{2equiv}
Let $f,g\in M(X,Z)$. With the notations above:
$$
f^{\ast }(\tilde {\alpha})=g^{\ast }(\tilde {\alpha})
\text { if and only if }
f^{\ast }(\alpha)=g^{\ast }(\alpha)
$$
and similarly for $\tilde {\beta}$.
\end{lemma}
\begin{proof}
One implication is obvious. In the opposite direction, we have
$$f^{\ast}(\alpha)=g^{\ast}(\alpha)\text{ \, and \, }
f^{\ast}(\overline {\alpha})=g^{\ast}(\overline {\alpha}).$$
Therefore $f^{\ast}(\tilde {\alpha})-g^{\ast}(\tilde {\alpha})$ is
a $(1,1)$ integral element, so it does not belong to the
transcendental lattice.
\end{proof}

Let us consider the injection
$$
\begin{aligned}
M_r(X,Z)/\sim &\hookrightarrow (T_{X,\mathbb Z}\times T_{X,\mathbb Z} )\otimes \mathbb {R} \\
[f] &\longmapsto v_f:=(\frac 1{\|\tilde {\alpha } \|}f^{\ast}(\tilde {\alpha }),
\frac 1{\|\tilde {\beta}\|}f^{\ast}(\tilde {\beta })).
\end{aligned}
$$

Then:
\begin{proposition} \label{dist2}
Let $f,g\in M_r(X,Z)$ such that $g\notin [f]$. Then:
$$
\|v_f-v_g\|\ge \frac 1{2\sqrt r}
$$
\end{proposition}
\begin{proof}
Assume first that $f^{\ast }(\tilde {\alpha })\ne g^{\ast}(\tilde
{\alpha })$. Then, by arguing as in \ref{dist} we obtain
$$
\|v_f-v_g\|\ge \|\frac 1{\|\tilde {\alpha }\|}f^{\ast}(\tilde
{\alpha })- \frac 1{\|\tilde {\alpha }\|}g^{\ast}(\tilde {\alpha
})\|\ge \frac 1{\sqrt {r}}> \frac 1{2\sqrt {r}}.
$$

If $f^{\ast }(\tilde {\alpha } )= g^{\ast}(\tilde {\alpha })$, then $f$ and $g$ coincide on
$\left\langle \tilde {\alpha } \right\rangle $, which implies
$\left\langle \tilde {\alpha } \right\rangle \ne T_X$.

Observe that
$$
((f_{\ast }-g_{\ast })(f^{\ast }(\tilde {\beta })-g^{\ast }(
\tilde {\beta } )),\tilde {\alpha })=(f^{\ast }(\tilde {\beta }
)-g^{\ast }({\tilde \beta }), f^{\ast }(\tilde {\alpha } )-g^{\ast
}(\tilde {\alpha } ))=0
$$
and
$$
((f_{\ast }-g_{\ast })(f^{\ast }(\tilde {\beta} )-g^{\ast }(\tilde {\beta })),\tilde {\beta})=
\|f^{\ast }(\tilde {\beta })-g^{\ast }(\tilde {\beta })\|^2\ne 0.
$$

Hence, $(f_{\ast }-g_{\ast })(f^{\ast }(\tilde {\beta })-g^{\ast
}(\tilde {\beta } ))$ is a non trivial element in the lattice
$R_{\mathbb Z}$, being $R=\left\langle \tilde {\alpha
}\right\rangle ^{\perp }$ the orthogonal Hodge structure to $
\left\langle \tilde {\alpha } \right\rangle $ in $T_X$. Hence its
norm is greater or equal to $\|\tilde {\beta }\|$. We get
$$
\begin{aligned}
\|\tilde {\beta }\| &\le \|(f_{\ast }-g_{\ast })(f^{\ast }(\tilde {\beta })-g^{\ast }(\tilde {\beta } ))\| \\
&\le \| f_{\ast }(f^{\ast }(\tilde {\beta })-g^{\ast }(\tilde {\beta } ))  \|+
\|g_{\ast}(f^{\ast }(\tilde {\beta })-g^{\ast }(\tilde {\beta } ))\|\le
2\sqrt{r}\|f^{\ast }(\tilde {\beta })-g^{\ast }(\tilde {\beta } )\|
\end{aligned}
$$
and the proposition follows.
\end{proof}

Finally, by using the packing lemma with $R=\sqrt {2r}$, $d=\frac
1{4\sqrt{r}}$ and the fact that $r\le \rho$ (see \ref{degrat}) and
\ref{tgl2} we have:
$$
\begin{aligned}
m(X,Z)\le & 4K_Z^2\,K_X^2 \,\sum_{r=1}^{\rho } \# (M_r(X,Z)/\sim)  \\
      \le & 4K_Z^2\,K_X^2 \,\rho \, P(4\sqrt 2\, \rho,2\dim T_X)
      \\
      \le &4K_X^2\,K_X^2\, P(4\sqrt 2\, \rho,2b_2(X)-2),
\end{aligned}
$$
the last inequality comes from $\dim T_X\le b_2(X)-1$. Therefore the proof of \ref{th2} is finished.

\section{Threefolds}

We now consider the $3$-dimensional case.   As we will see below,
we can concentrate on the geometric part of the proof, since the
representation of $M_r(X,Z)/\sim $ in the transcendental lattice
and the estimation of the distance work, word by word, in the same
way.

We assume $X,Z$ of general type, $K_X,\, K_Z$  are nef, $p_g(Z)\ge 2$ and
$$
\dim(\varphi_{|2K_Z|}(Z))\ge 2.
$$
Fix two linearly independent
$(3,0)$ forms $\alpha $ and $\beta $. As in the precedent
sections, given $f: X - - \to Z$ dominant,  we focus in the
estimation of the number of elements of
$$
[f]=\{g:X - - \to Z \,|\, f^{\ast }(\alpha )=g^{\ast }(\alpha ),\,
f^{\ast }(\beta )=g^{\ast }(\beta)\}.
$$

\begin{remark}
We use a pencil on $Z$ to reduce the proof to the case of
surfaces. We could instead fix $3$ forms and try to reduce
directly to curves. This method fails, since the corresponding map
to $\mathbb P^2$ could not be dominant. Observe that we cannot
choose generic forms since in order to apply packing arguments we
need to fix them with some minimal properties.
\end{remark}

We follow closely the case of surfaces: we have a pencil on $Z$,
$\beta $ is a general element of the pencil and its divisor of
zeros is:
$$
B+S_1+\dots +S_s,
$$
where $B$ is the base divisor.

 In the same way, the divisor of
zeros of $f^{\ast }(\beta )$ can be written:
$$
B'+S'_1+\dots +S'_{s'},
$$
where $B'$ is the base divisor.

Denote $r=deg(f)$. Then:

\begin{lemma}
One has the following inequality:
$$
s\le K_Z^3.
$$
\end{lemma}
\begin{proof}
 Since $S_i$ moves, a
convenient pluricanonical map sends $S_i$ to a surface. Therefore,
for $l>>0$, $(lK_Z)^2S_i>0$, so $K_Z^2S_i>0$. Hence, by the
nefness of $K_Z$:
$$
s\le \sum _{i=1}^{s}K_Z^2S_i = K_Z^2(K_Z-B)\le K_Z^3.
$$
\end{proof}

Consider $Z' \longrightarrow Z$ to be the minimal resolution of
the pencil $Z - - \to \mathbb P^1$ induced by $\alpha $ and $\beta
$ and let $X'$ be the minimal resolution of the induced pencil on
$X\times _Z Z'$. The general fibre of the pencil on $Z'$ is a
disjoint union of smooth surfaces $T_1,\dots ,T_s$, being $T_i$ a
desingularization of $S_i$. We have in the same way the smooth
surfaces $T'_1,\dots , T'_{s'}$ on $X'$. Then the map $f$ and the
forms $\alpha $, $\beta $ pullback to $f'$, $\alpha '$, $\beta '$
and we obtain
$$[f]=[f']=\{g':X' - - \to Z' \,|\, f'^{\ast }(\alpha ')=g'^{\ast
}(\alpha '),f'^{\ast }(\beta ')=g'^{\ast }(\beta ')\}.$$

Now we divide the set $[f']$ into subsets depending on the image
of the fixed surface $T_1'$:
$$
[f']=\bigcup _i\{g\in [f']\,|\,g(T_1')=T_i\}\subset \bigcup _i
M(T_1',T_i).
$$
The second inclusion follows since the surface $T_1'$ intersects
the rational domain for $X -- \to Z$.

\begin{proposition} \label{tgl3}
One has:
$$
\#[f]\le  4\cdot 9^2\, K_Z^3\,h^2\, P(36\sqrt 2 \,h,20 h +6),
$$
where $h=h^0(X,\mathcal O(2K_X))+h^0(X,\Omega ^2_X)-p_g(X)+1$.
\end{proposition}
\begin{proof}
By the inclusion above
$$
\#[f]=\#[f']\le \sum_{i=1}^{s}m(T_1',T_i).
$$
The surfaces $T'_1,T_i$ are of general type, since they move in a
rational pencil and the threefolds are of general type.

Observe that, since the image of $\varphi _{|2K_Z|}$ has dimension
$\ge 2$ , there exist at least two elements $\alpha_1,\alpha_2\in
H^0(Z',\omega_{Z'}^{\otimes 2})$ such that the residues
$$
\text{Res }_{T_i}(\frac{\alpha _1}{\beta '}),
\qquad \text{Res }_{T_i}(\frac{\alpha _2}{\beta '})
$$
define on each component $T_i$ two linearly independent
holomorphic $2$-forms. Therefore $p_g(T_i)\ge 2$. With these
hypothesis we can apply corollary \ref{cor2} to obtain
$$
m(T_1',T_i)\le 4\cdot 9^2\chi ^2\,P(36\sqrt 2 \,\chi ,20 \chi
+6),
$$
where $\chi $ is $\chi  (\mathcal O_{T_1'})$.

 To finish the proof
we have to bound $\chi $ by $h$ and use $s\le K_Z^3$.

Let us consider the exact sequence of sheaves on $X'$:
$$
0\to \omega _X'  \longrightarrow \omega _X' (T_1') \longrightarrow
\omega _{T_1'} \to 0.
$$
By taking the attached long exact sequence in cohomology we obtain
$$
p_g(T_1')\le h^1(X,\omega _X)+h^0(X,\omega ^{\otimes 2})-p_g(X) =
h^0(X,\Omega ^1_{T_1'})+h^0(X,\omega ^{\otimes 2})-p_g(X).
$$
Since $\chi \le p_g(T'_1)+1$ we are done.
\end{proof}

To prove \ref{thn} we can imitate the proof of  \ref{th2} given in
the last section. The only difference is that the analogous
statement to Proposition \ref{2equiv} is no longer true. However
the obvious implication
$$
f^{\ast }(\tilde {\alpha })=g^{\ast }(\tilde {\alpha })
\,\Rightarrow \, f^{\ast } (\alpha )=g^{\ast }(\alpha )
$$
is enough to ensure that the equivalence class of $f$ is contained
in $[f]$.

Then, by using \ref{tgl3}:
$$
m(X,Z) \le  4\cdot 9^2\, K_Z^3\,h^2\, P(36\sqrt 2 \,h,20 h +6)\, \rho \, P(4\sqrt 2 \rho ,2\text {dim }T_X).
$$
The statement of  \ref{th3} follows replacing $\rho$ with $\frac{K_X^3}{K_Z^3}$.

\section{Torsion lemma}

In this section we generalize the torsion part of Tanabe's work to
higher dimensions. With some hypotheses, this allows us to produce
bounds for $m(X,Z)$ in any dimension.

Let $f:Z\to X$ and $g:Z\to X$ dominant map of degree $r.$ We let
$f_{\ast}$, $f^{\ast},$ $f_{d},$ $g_{\ast},$ $g^{\ast}$ and
$g_{d}$ induced maps (see section 2).

 We have the following:

\begin{lemma} \label{torsione}
If $f_{d}=g_{d}$ for some $d>2r$ then $f^{\ast}=g^{\ast}$.
\end{lemma}
\begin{proof}
Let $h=f^{\ast}-g^{\ast}$ we have to prove that  $T_Z=\ker (h).$
If not, let $V$ be Hodge polarized structure orthogonal to $\ker
h.$ let $h':V \to T_{Z}$ be the  restriction of $h.$ Now $h'$ is
injective. Set $\mu \in V_{\mathbb Z}$ such that its norm is
minimal in the lattice. We consider
$$\lambda= h'(\mu )= f^{\ast}(\mu ) -g^{\ast}(\mu ).$$
We would have that  $\lambda\neq 0.$  Moreover from the hypothesis
$f_{d}=g_{d}$ we have that $\lambda= d\cdot\sigma $ where $\sigma
\in T_{\mathbb Z,Z}$ is an integral class.

We also consider $\beta_{f}=f_{\ast}(\lambda) $ and $\beta_{g}=
g_{\ast}(\lambda).$ We have that $\beta_{f}$ (and $\beta_{g)}$)
are in $V.$ To see this, first we remark that
$f_{\ast}f^{\ast}=g_{\ast}g^{\ast},$ since
 $$
 f_{\ast}f^{\ast}(\alpha)=r\cdot \alpha=g_{\ast}g^{\ast}(\alpha).
 $$
Now fix $\alpha \in \ker(h),$ $f^{\ast}(\alpha)=g^{\ast}(\alpha),$
we have
$$
\begin{aligned}
(\alpha,& \beta_{f})=(\alpha, f_{\ast}(\lambda))=
(f^{\ast}(\alpha),\lambda)=(f^{\ast}(\alpha), f^{\ast}(\mu)
-g^{\ast}(\mu))= \\
 &= (f^{\ast}(\alpha),f^{\ast}(\mu))-
 (f^{\ast}(\mu ),g^{\ast}(\mu ))=(f_{\ast}f^{\ast}(\alpha),\mu ))-
(f^{\ast}(\alpha),g^{\ast}(\mu ))= \\
&= (g_{\ast}g^{\ast}(\alpha),\mu ))-
(f^{\ast}(\alpha),g^{\ast}(\mu ))=(g^{\ast}(\alpha),g^{\ast}(\mu
))- (f^{\ast}(\alpha),g^{\ast}(\mu ))\\
 &= (h(\alpha),g^{\ast}(\mu
)) =0.
\end{aligned}
$$

Then we have
 $$\beta_{f}-\beta_{g}= (f_{\ast}-g_{\ast})(\lambda)
 = (f_{\ast}-g_{\ast})(f^{\ast} -g^{\ast})(\mu )$$
is not zero. Indeed
$$
(\beta_{f}-\beta_{g}, \mu )= ((f^{\ast}-g^{\ast})(\mu
),f^{\ast}-g^{\ast})(\mu ))= \|\lambda \|^2\neq 0.
$$
 It follows that either $\beta_f$ or
$\beta_{g}$ are not zero.

We may assume now that  $\beta_f=f_{\ast}(\lambda)\neq 0.$ Recall
that we have that $\lambda=d\cdot \sigma$ where $ \sigma \in
T_{\mathbb Z,Z}.$ We have then
$$
\|f_{\ast}((f^{\ast}-g^{\ast})(\mu))\|= \|f_{\ast}(\lambda)\|=
d\cdot\|f_{\ast}(\sigma)\|\geq  d\cdot\|\mu\|,$$ by the minimality
of $\|\mu \|$.

In addition:
$$
\|f_{\ast}((f^{\ast}-g^{\ast})(\mu))\|\leq \sqrt
r\,\|(f^{\ast}-g^{\ast})(\mu))\|\leq \sqrt r\, (\|f^{\ast}(\mu)\|+
\|g^{\ast}(\mu)\|)=
 2r \|\mu\|. $$
Hence $$ d\leq 2r.$$
 \end{proof}

The rest of the section is devoted to the proof of Theorem
\ref{thn}. We fix $X,Z$ two $n$-dimensional varieties of general
type, $n\ge 2$, such that $K_Z$ is nef and $\text{dim }(\varphi
_{|K_Z|}(Z))\ge n-1$ (in particular $p_g(Z)\ge n$).

\begin{definition}
We say that two maps $f,g\in M(X,Z)$ are equivalent if $f^{\ast }=g^{\ast }$ on $T_{Z}$.
\end{definition}

As usual we would like to compute the number of elements of the
class $[f]$ of a map $f$. We consider a general projection of the
image of the canonical map of $Z$. Then we obtain a rational
dominant map $\phi :Z - - \to \mathbb P^{n-1}$. By definition
$\phi \circ f=\phi \circ g$.

Observe $\phi $ can be written as $Z - - \to \mathbb P(V^{\ast
})$, where $V$ is a $n$-dimensional vector space contained in
$H^0(Z,\omega_Z(-F))$, $F$ being the fixed divisor of the linear
system attached to $V$.
 The general fibre  of $\phi $ is
$$
C_1+\dots +C_s
$$
and can be thought of as the common zeros of $\{s_1,\dots
,s_{n-1}\},$ where $s_i \in V$. Let $s_n$ be another element in
$V$ such that $s_1,\dots ,s_n$ is a basis of $V$.

The fibre of $\phi \circ f$ is of the form
$$
C_1'+\dots +C'_{s'}.
$$

We consider a resolution $\pi :Z' \longrightarrow Z$ of the singularities of $\phi$.
We put $\pi ^{\ast }(s_i)=s_i'\cdot s_{E_0}$, where $E_0$ is the fixed divisor of the pull-back of the
linear system and $s_{E_0}$ is an equation for this divisor. Then
$$
\left\langle s_1',\dots ,s_n'\right\rangle \subset
H^0(Z',\pi^{\ast }(\omega _Z(-F))(-E_0))
$$
defines the map $\phi'=\phi \circ \pi: Z'\longrightarrow \mathbb
P^n$. The normalizations $N_i$ of $C_i$ are the components of the
general fibres of $\phi '$.  By construction we see
$N_1+\dots+N_s$ as the common zeros of $\{s_1',\dots,s_{n-1}'\}$.
In particular we have the rational equivalence of $1$-cycles,
$$
N_1+\dots +N_s \sim_{rat} (\pi^{\ast }(K_Z)-E)^{n-1},
$$
where $E=E_0+\pi ^{\ast}(F)$ is an effective divisor (and $h^0(\pi^{\ast }(K_Z)-E)>0$ by construction).

Notice also that $s_n'^{\otimes n}$ restricts to $N_i$ giving a holomorphic form $\alpha $.
Locally this  form is computed as the following residue:
$$
\text{ Res }_{N_i}(\frac {s'_n\cdot \dots \cdot s'_n \cdot s_{E_0}}{s_1'\cdot \dots \cdot s'_{n-1}}).
$$
Analogously we can resolve the singularities of the map $X - - \to \mathbb P^n$ and the
general fibre is $N_1'+\dots +N_{s'}'$, being $N_i'$ the desingularization of $C_i'$.

Then, since $N_1'$ intersects the rational domain of $X$ and $Z$
(defined in \ref{ratdom}):
$$
[f]=\bigcup _{i=1}^s \{h:N_1'\longrightarrow N_i\,|\, h^{\ast }(\alpha) \text { fixed}\}.
$$
By \ref{tgl1} we obtain
$$
\#[f]\le s\cdot 4(g(N_1')-1).
$$
To go further, we need the following
\begin{lemma}
With our hypothesis
\begin{enumerate}
\item[a)] $s\le K_Z^n.$
\item[b)] $g(N_i')\le \frac {nK_X^n}2+1$,  $g(N_i)\le \frac {nK_Z^n}2+1.$
\end{enumerate}
\end{lemma}
\begin{proof}
\begin{enumerate}
\item[a)] Since  $C_i$ moves, $K_Z\,C_i\ge 1$. Then
$$
s\le K_Z\, (\sum _{i=1}^s C_i)=\pi ^{\ast }(K_Z)(\sum _{i=1}^s N_i)=\pi^{\ast}(K_Z)(\pi^{\ast }(K_Z)-E)^{n-1}.
$$
To see $s\le K^n_Z$, it is enough to prove that
$$
\pi^{\ast}(K_Z)(\pi^{\ast }(K_Z)-E)^{n-1}\le K^n_Z,
$$
which follows from the positivity of
$$
\begin{aligned}
&\pi ^{\ast }(K_Z)(\pi ^{\ast }(K_Z)^{n-1}-(\pi ^{\ast }(K_Z)-E)^{n-1})= \\
& \pi ^{\ast }(K_Z)E(\pi ^{\ast }(K_Z)^{n-2} +\pi ^{\ast
}(K_Z)^{n-3}(\pi ^{\ast }(K_Z)-E)+\dots )
\end{aligned}
$$
(using the fact  that $K_Z$ is nef).
 \item[b)] We give the proof
on $Z$. Observe that $p_a(C_i)\le p_a(\sum _i C_i)$ since
$p_a(C_1)=\dots =p_a(C_s)\ge 2$ and all the curves are reduced.
Therefore
$$
g(N_i)\le p_a(C_i) \le p_a(\sum _i C_i)
$$
and
$$
\begin{aligned}
2p_a(\sum _i C_i)-2 =&(K_Z+(n-1)(K_Z-F))(K_Z-F)^{n-1} \\
=&(nK_Z-(n-1)F)(K_Z-F)^{n-1}\le n K_Z^n.
\end{aligned}
$$
The last inequality is proved as in the first part.
We are done.
\end{enumerate}
\end{proof}

A direct consequence of the lemma and the discussion above is the inequality:
\begin{proposition}\label{tgln}
We have the bound:
$$
\#[f]\le 2nK_X^nK_Z^n
$$
\end{proposition}

Now we fix a degree $r$ and we consider the map
$$
M_r(X,Z)/\sim \longrightarrow Hom (T_{Z,\mathbb Z}/(2r+1)T_{Z,\mathbb Z},T_{X,\mathbb Z}/(2r+1)T_{X,\mathbb Z})
$$
which sends $f$ to $f_d$. The injectivity is an application of \ref{torsione}.
Then, by \ref{tgln}:
$$
\begin{aligned}
m(X,Z)&\le  2n \,K_X^n K_Z^n\, \sum _{i=1}^{\rho }(2r+1)^{\text{dim}T_Z \cdot \text{dim}T_X}
\\
&\le 2n\, K_X^n K_Z^n \,
\rho \,(2\rho +1)^{b_n(X)\, b_n(Z)}.
\end{aligned}
$$
This finishes the proof of
\ref {thn}.

\vskip 10mm

\bibliographystyle{amsplain}

\newcommand{\noopsort}[1]{} \newcommand{\printfirst}[2]{#1}
  \newcommand{\singleletter}[1]{#1} \newcommand{\switchargs}[2]{#2#1}
\providecommand{\bysame}{\leavevmode\hbox to3em{\hrulefill}\thinspace}
\providecommand{\MR}{\relax\ifhmode\unskip\space\fi MR }
\providecommand{\MRhref}[2]{%
  \href{http://www.ams.org/mathscinet-getitem?mr=#1}{#2}
}
\providecommand{\href}[2]{#2}

\vskip 5mm
\noindent
{\sc Juan Carlos Naranjo} \\
Departament d'\`Algebra i Geometria, Universitat de Barcelona \\
Gran Via 585, 08007 Barcelona, Spain \\
{\tt jcnaranjo@ub.edu}

\vskip 3mm
\noindent {\sc Gian Pietro Pirola}\\
Dipartimento di Matematica, Universit\`a di Pavia\\
via Ferrata 1, 27100 Pavia, Italia\\
{\tt pirola@dimat.unipv.it}

\end{document}